\documentclass[12pt]{amsart}
\usepackage{amsmath}

\def\p{\partial}

\newcommand\ce{{\mathbb C}}

\def\G{g_{ij}}

\def\K{K\"ahler }
\def\KR{K\"ahler-Ricci }

\def\KRS{K\"ahler-Ricci soliton }
\def\KRS{K\"ahler-Ricci soliton }

\def\be{\begin{equation}}
\def\ee{\end{equation}}

\def\p{\partial}

\newtheorem{thm}{Theorem}

\theoremstyle{definition}

\theoremstyle{remark}
\newtheorem{rem}{Remark}
\numberwithin{equation}{section}

\begin{document}

\title{A note on the uniformization of gradient K\"ahler  Ricci solitons}

%    Information for first author
\author{Albert Chau$^{1}$}
%    Address of record for the research reported here
\address{Harvard University, Department of Mathematics,
  One Oxford Street, Cambridge, MA 02138, USA}
\email{chau@math.harvard.edu}
%    Current address
%\curraddr{}
%    \thanks will become a 1st page footnote.
%\thanks{}

%    Information for second author
\author{Luen-Fai Tam$^2$}
\thanks{$^{1}$Research
partially supported by The Institute of Mathematical Sciences, The Chinese
University of Hong Kong, Shatin, Hong Kong, China.}
\thanks{$^2$Research
partially supported by Earmarked Grant of Hong Kong \#CUHK4032/02P}

\address{Department of Mathematics, The Chinese University of Hong Kong,
Shatin, Hong Kong, China.} \email{lftam@math.cuhk.edu.hk}

%    General info
\renewcommand{\subjclassname}{%
  \textup{2000} Mathematics Subject Classification}
% Global differential geometry
% 53C21 Methods of Riemannian geometry, including PDE methods;
%   curvature restrictions [See also 58J60]
% 53C44 Geometric evolution equations (mean curvature flow)
% 53C55 Hermitian and K\"ahlerian manifolds [See also 32Cxx]
% Qualitative properties of solutions
% 35B35 Stability, boundedness
% Partial differential equations on manifolds; differential operators
% 58J37 Perturbations; asymptotics
% Systems theory; control - Stability
% 93D05 Lyapunov and other classical stabilities
%       (Lagrange, Poisson, $L^p, l^p$, etc.)

\date{April  2004 (Revised in July 2004).}

%\dedicatory{}

%\keywords{}

\begin{abstract}
Applying a well known result for attracting fixed points of
biholomorphisms \cite{RR, V}, we observe that one immediately obtains the following result: if $(M^n,g)$ is a complete non-compact gradient K\"ahler-Ricci soliton which is either steady with positive Ricci curvature so that the scalar curvature attains its maximum at some point, or expanding with
non-negative Ricci curvature, then $M$ is biholomorphic to $\ce ^n$.
\end{abstract}

\maketitle \markboth{Albert Chau and Luen-Fai Tam} {A note on Gradient
\K Ricci solitons.}

We will show the following:

\begin{thm}\label{t}
If $(M^n,g)$ is a complete non-compact gradient K\"ahler-Ricci soliton which is either steady with positive Ricci curvature so that the scalar curvature attains its maximum at some point, or expanding with
non-negative Ricci curvature, then $M$ is biholomorphic to $\ce ^n$.
\end{thm}

Recall that a K\"ahler manifold
$(M,g_{i \bar j}(x))$   is said to be   a \KR soliton if there is a family of
biholomorphisms $\phi_t$ on $M$, given by a holomorphic vector field $V$, such that
$\G(x,t)=\phi_t^*(\G(x))$ is  a solution of the \KR flow:
\begin{equation}\label{1.e1}
\begin{split}
\frac{\p}{\p t}g_{i \bar j}&=-R_{i\bar j}-2\rho g_{i \bar j}\\
g_{i \bar j}(x,0)&=g_{i \bar j}(x)
\end{split}
\end{equation}
for $0\le t<\infty$, where
$R_{i\bar j}$ denotes the Ricci tensor at time $t$ and $\rho$ is a constant. If
$\rho=0$, then the \KR soliton  is said to be of {\it steady type} and if $\rho>0$ then
the \KR soliton  is said to be of {\it expanding type}. We always
assume 
that $g$ is
complete and $M$ is non-compact.
If in addition, the holomorphic vector field is
given by  the gradient of a real valued function $f$, then it is called a
gradient \KR soliton. Note that in this case, we have that
\begin{equation}\label{1.e4}
\begin{split}
f_{i\bar j}&= R_{i\bar j}+2\rho  g_{i\bar j}\\
f_{ij}&=0.
\end{split}
\end{equation}
If $(M,g)$ is a gradient \KR soliton (of steady or expanding type) which is either steady with positive Ricci curvature so that the scalar curvature attains its maximum at some point, or expanding with
non-negative Ricci curvature, then one can show that $\phi_t$, the flow
on $M$ along the vector field  $\nabla f$, satisfies: 
   
\begin{enumerate}
\item[(i)] $\phi_t$ is a biholomorphism from $M$ to $M$  for all $t \geq 0$,
\item[(ii)]  $\phi_t$ has a unique fixed point $p$, i.e. $\phi_t(p)=p$
  for all $t \geq 0$,
\item[(iii)] M is attracted to $p$ under $\phi_t$ in the sense that
  for any open neighborhood $U$ of $p$ and for any compact subset $W$ of $M$, there exists $T>0$ such that
$\phi_t(W)\subset U$  for all $t \geq T$.\\
\end{enumerate}

Condition (i) is clear.  Condition (ii) is shown in \cite{CT, Ha4}.
To see that condition (iii) holds, we consider any $R>0$ and let $B(R)$ be the
geodesic ball of radius $R$ with center at $p$ with respect to the metric $g(0)$.
From the proof of Lemma 3.2 in \cite{CT}, there exists $C_R>0$ such that for any
$q\in B(R)$ and for any $v\in T^{1,0}(M)$ at $q$,
$$
||v||_{\phi_t^*(g)}\le \exp(-C_Rt)||v||_g.
$$
Since $\phi_t(p)=p$, it is easy to see that given any open set  $U\subset M$
containing $p$, we have $\phi_t(B(R))\subset U$ provided $t$ is large,
and thus condition (iii) is satisfied.

The following theorem was proved for the case $M=\ce^n$
in \cite{RR}, and was later observed to be true on a general complex
manifold $M$ in \cite{V}.    
\begin{thm}\label{v}
Let F be a biholomorphism from a complex manifold $M^n$ to itself and let $p
\in M^n$ be a fixed point for F.  Fix a complete Riemannian metric g
on M and define $$\Omega:=\{x\in M: \lim_{k\to \infty} dist_g (F^k(x), p)=0\}$$where  $F^k=F\circ F^{k-1}, F^1=F$.
Then  $\Omega$ is biholomorphic to $\ce^n$ provided $\Omega$ contains an
 open neighborhood around p. 
\end{thm}

\begin{proof}[Proof of Theorem 1]
By conditions (i)-(iii) we may apply Theorem \ref{v} to
the biholomorphism $\phi_1:M \to M$ to conclude
that $M$ is biholomorphic to $\ce^n$. 
\end{proof}

\begin{rem}In the first version of this article we proved Theorem \ref{v} in a special case.  We would like to thank Dror Varolin
 for pointing out to us that what we proved had been known earlier \cite{RR, V}.
\end{rem}  

\begin{rem}
After posting the first version of this article we learned that
Theorem \ref{t} in the case of a steady gradient \K Ricci soliton had been  known independently to Robert Bryant \cite{B}. 
\end{rem}

\bibliographystyle{amsplain}

\begin{thebibliography}{10}
\bibitem{B} Bryant, Robert, {\sl Gradient Kahler Ricci Solitons},
  arXiv eprint 2004. arXi:math.DG/0407453 
\bibitem{CT} Chau, Albert and Tam, Luen-fai, {\sl Gradient \KRS and a 
uniformization conjecture}, arXiv eprint 2003. arXi:math.DG/0310198.

\bibitem{Ha4} Hamilton, Richard S.,
 {\sl Formation of Singularities in the Ricci Flow},

  Contemporary Mathematics,  \textbf{71} (1988), 237-261.


\bibitem{RR}Rosay, J.P. and Rudin, W.,{\sl Holomorphic Maps from $\ce^n$
  to $\ce^n$}, Trans. AMS 310 (1988), 47-86.

\bibitem{Stern} Sternberg, Shlomo, {\sl Local contractions and a theorem  
of Poincaré},   Amer. J. Math., \textbf{79}  (1957), 809--824.

\bibitem{V} Varolin, Dror, {\sl The density property for complex
  manifolds and geometric structures II}, Internat. J. Math. \textbf{11} (2000), no. 6, 837--847. 
 



\end{thebibliography}

\end{document}